\documentclass[12pt,a4paper,reqno]{article}    % Specifies the document style.
\usepackage{graphics}
\usepackage{cite}
\usepackage[demo]{graphicx}
\usepackage[fleqn]{amsmath}
\usepackage{booktabs}
\usepackage{multirow}
\usepackage{epsfig}
\usepackage{epstopdf}
\usepackage{booktabs}
\usepackage{placeins}
\usepackage{anysize}
\usepackage{enumerate}
\usepackage{latexsym,amssymb,amsmath,amscd,amsthm,amsxtra}
\usepackage{subfigure}
\usepackage{tabularx}
\usepackage[table]{xcolor}
\allowdisplaybreaks
\definecolor{lightgray}{gray}{0.9}
\newtheorem{theorem}{Theorem}
\newtheorem{prop}{Proposition}
\newtheorem{exmp}{Example}

\numberwithin{equation}{section}

\def\({\left( }
\def\){\right )}
\begin{document}
\title {\textbf{A Finite Difference Scheme based on Cubic Trigonometric B-splines for Time Fractional Diffusion-wave Equation }}
\author{Muhammad Abbas$^1$\thanks{Corresponding authors: m.abbas@uos.edu.pk}\\
{$^1$\small Department of Mathematics, University of Sargodha, Sargodha, Pakistan.}\\
}
\date{}
\maketitle
\begin{abstract}
In this paper, we propose an efficient numerical scheme for the approximate solution of the time fractional diffusion-wave equation with reaction term based on cubic trigonometric basis functions. The time fractional derivative is approximated by the usual finite difference formulation and the derivative in space is discretized using cubic trigonometric B-spline functions. A stability analysis of the scheme is conducted to confirm that the scheme does not amplify errors. Computational  experiments are also performed  to further establish the accuracy and validity of the proposed scheme. The results obtained are compared with a finite difference schemes based on the Hermite formula and radial basis functions. It is found that our numerical approach performs superior to the existing methods due to its simple implementation, straight forward interpolation and very less computational cost.  
\end{abstract}
\maketitle
%%%%%%%%%%%%%%%%%%%%%%%%%%%%%%%%%%%%%%%%%%%%%%%%%%%%%%%%%%%%%%%%%
% Type in your PAPER, starting below:
\section{Introduction}
\subsection{Problem Description}
For $T>0$ and $\Omega=[a,b]$,  we consider the following model of the time fractional diffusion-wave equation with reaction term,
\begin{equation}\label{1.1}
    \frac{\partial^ \gamma }{\partial t^\gamma} u(x,t)+\alpha u(x,t)= \frac{\partial ^2}{\partial x^2} u(x,t)+f(x,t)\,\,\,\,\,\,\,\,\,\,\,1 < \gamma \leq 2,\,\,\,\,x \in \Omega,\,\,\,0\leq t\leq T
\end{equation}
with initial conditions
 \begin{equation}\label{1.2}
\left\{
\begin{array}{cc}
 u(x,0)=\phi_1(x),\\
 u_t(x,0)=\phi_2(x)
\end{array}\qquad\qquad\qquad x \in \Omega, \right.
\end{equation}

and the following boundary conditions
\begin{equation}\label{1.3}
\left\{
\begin{array}{cc}
 u(a,t)=\psi_1(t),\\
 u(b,t)=\psi_2(t)
\end{array}\qquad\qquad\qquad 0\leq t\leq T,\right.
\end{equation}
where $a,b, \phi_1(x),\phi_2(x),\psi_1(t)$ and $\psi_2(t)$ are given, $\alpha>0$ is the reaction coefficient and  $\frac{\partial^\gamma}{\partial t^\gamma}u(x,t)$ represents the Caputo fractional derivative of order $\gamma$  given by \cite{pod}
\begin{equation}\label{1.4}
\frac{\partial^\gamma}{\partial t^\gamma}u(x,t)=
\begin{cases}
\frac{1}{\Gamma(2-\gamma)}\int\limits_0^t \frac{\partial^2 u(x,s)}{\partial^2 s}(t-s)^{1-\gamma}ds,\qquad 1<\gamma <2\\
\frac{\partial^2 u(x,s)}{\partial^2 s},\qquad\qquad\qquad\qquad\qquad\qquad \gamma=2.
\end{cases}
\end{equation}
To obtain the time fractional diffusion-wave equation from the standard diffusion or wave equation we replace the ordinary first or second time derivative by a fractional derivative of order $\gamma$ where $0<\gamma <1$ or $1<\gamma <2$. As $\gamma$ changes from $0$ to $2$ the process transforms from slow diffusion to classical diffusion and diffusion-wave to classical wave phenomenon. We consider, in this paper, the case of diffusion-wave i.e. $1 < \gamma < 2$. It can be used to deal with viscoelastic problems and disordered media to examine structures, semiconductors and dielectrics.
\subsection{Applications and Literature Review}
The subject of fractional calculus \cite{pod,mai,hil,kilbas} in its modern form has a history of at least three decades and has developed rapidly due to its wide range of applications in  fluid mechanics, plasma physics, biology, chemistry, mechanics of material science  and so on \cite{hil,sok}. Other applications include   system control \cite{pod}, viscoelastic flow \cite{die},  hydrology \cite{bek,meer}, tumor development \cite{iom} and finance \cite{scale1,scale2,scale3}. Since the fractional models in certain situations tend to behave more appropriately than the conventional integer order models, several techniques have been developed to study these models. These techniques have been continuously improved and modified to achieve more and more accuracy.\\

Since exact analytical solutions of only a few fractional differential equations exist, the search for approximate solutions  is a concern of many recently published articles. Many research publications have been devoted to numerical techniques for solving time fractional diffusion-wave equations. Zeng \cite{zeng} proposed two second order stable and one conditionally stable finite difference schemes for the time fractional diffusion-wave model. Using a class of finite difference methods based on the Hermite formula,  Khader and Adel \cite{khader} has obtained numerical solutions of fractional diffusion-wave equation.  Avazzadeh \textit{et al} \cite{vahid} have obtained numerical solutions of fractional diffusion-wave equation by using Radial basis function method. Pskhu \cite{psku} has obtained fundamental solutions of fractional order diffusion-wave equation. It has been shown that this fundamental solution gives the corresponding solutions for diffusion and wave equations when the fractional order is equal to one or approaches  two. Povstenko \cite{pov} has discussed Neumann boundary-value problems for a time-fractional diffusion-wave equation in a half-plane.  Numerical solutions to the fractional diffusion-wave equation under Dirichlet and Neumann boundary conditions has been obtained by Povstenko.  Liemert and  Kienle \cite{lie} have discussed time fractional wave-diffusion equation in an inhomogeneous half-space.  Ren and Sun \cite{ren} have obtained efficient numerical solutions of the multi-term time fractional diffusion-wave equation by using a compact finite difference scheme with fourth-order accuracy. Jin \textit{et al} \cite{jin} have utilized a Galerkin finite element method to find approximate solution  for a multi-term time-fractional diffusion equation.\\

An efficient numerical scheme based on trigonometric cubic B spline functions is presented in this paper to find the approximate solutions of the time fractional diffusion-wave equation with reaction term. First, we discretize the Caputo time fractional derivative by the usual finite difference formula and then use trigonometric cubic B-spline basis to approximate derivatives in space. Trigonometric cubic B-spline functions provide better accuracy than the usual finite difference schemes due to its minimal support and $C^2$ continuity. Numerical experiments are carried out and the obtained results of a special case are compared with those of \cite{khader}. The comparison shows that the presented scheme has an accuracy up to $10^{-11}$  whereas the scheme discussed in \cite{khader} has accuracy of $10^{-5}$. The scheme is shown to be unconditionally stable using a procedure similar to Von-Neumann stability analysis, whereas the scheme of \cite{khader} is conditionally stable. Numerical experiments confirm the validity and efficiency of the algorithm.\\

 The outline of this paper is as follows: In section 2, we give temporal  discretization using forward finite difference scheme of Eq. (\ref{1.1}). In section 3, we present the derivation of the scheme for the fractional diffusion-wave equation using the trigonometric cubic B-spline functions.  The stability analysis of the proposed scheme is given in section 4. Computational experiments are conducted to check the efficiency and validity of the scheme and the numerical results are reported in section 5. The last section is devoted to the concluding remarks of the study.
\section{Temporal Discretization}
To find time discretization of Eq. (\ref{1.1}) we discretize the Caputo  time fractional derivative $\frac{\partial^{\gamma}u(x,t)}{\partial t^{\gamma}}$  appearing in the equation using the usual finite difference method. Following the standard notations we let $t_n=n \Delta t, n=0,1,\cdots M$, where $\Delta t=\frac{T}{M}$ is the time step. First we approximate the second order differential operators using a forward finite difference method as follows
\begin{equation}\label{2.1}
    \frac{\partial ^2 u(x,s)}{\partial t^2}=\frac{u(x,t_{n+1})-2u(x,t_{n})+u(x,t_{n-1})}{(\Delta t)^2}+O(\Delta t)^2
\end{equation}
Using (\ref{2.1}), we can obtain an efficient approximation to the fractional derivative $\frac{\partial^{\gamma}u(x,t)}{\partial t^{1-\gamma}}$ as follows:
\begin{eqnarray}\label{2.2}
    \frac{\partial^{\gamma} u(x,t_{n+1})}{\partial t^{\gamma}}&=& \frac{1}{\Gamma(2-\gamma)} \int_0^{t_{n+1}} \frac{\partial^2 u(x,s)}{\partial^2 s} \frac{ds}{(t_{n+1}-s)^{1-\gamma}} \nonumber \\
    &=&\frac{1}{\Gamma(2-\gamma)}\sum\limits_{j=0}^n \int_{t_j}^{t_{j+1}} \frac{\partial^2 u(x,s)}{\partial^2 s} \frac{ds}{(t_{n+1}-s)^{1-\gamma}} \nonumber \\
    &=&\frac{1}{\Gamma(2-\gamma)}\sum\limits_{j=0}^n \frac{u(x,t_{j+1})-2u(x,t_{n})+u(x,t_{n-1})}{{\Delta t}^2} \int_{t_j}^{t_{j+1}}  \frac{ds}{(t_{n+1}-s)^{1-\gamma}}+e_{\Delta t}^{n+1} \nonumber \\
    &=&\frac{1}{\Gamma(2-\gamma)}\sum\limits_{j=0}^n \frac{u(x,t_{j+1})-2u(x,t_{n})+u(x,t_{n-1})}{{\Delta t}^2} \int_{t_{n-j}}^{t_{n+1-j}}  \frac{dr}{r^{1-\gamma}}+e_{\Delta t}^{n+1} \nonumber \\
    &=&\frac{1}{\Gamma(2-\gamma)}\sum\limits_{j=0}^n \frac{u(x,t_{j+1})-2u(x,t_{n})+u(x,t_{n-1})}{{\Delta t}^2} \int_{t_j}^{t_{j+1}}  \frac{dr}{r^{1-\gamma}}+e_{\Delta t}^{n+1} \nonumber \\
    &=&\frac{1}{\Gamma(3-\gamma)}\sum\limits_{j=0}^n \frac{u(x,t_{n+1-j})-2u(x,t_{n-j})+u(x,t_{n-1-j})}{\Delta t^\gamma} ((j+1)^{2-\gamma}-j^{2-\gamma})+e_{\Delta t}^{n+1} \nonumber \\
    &=&\frac{1}{\Gamma(3-\gamma)}\sum\limits_{j=0}^n b_j \frac{u(x,t_{n+1-j})-2u(x,t_{n-j})+u(x,t_{n-1-j})}{\Delta t^\gamma} +e_{\Delta t}^{n+1} \nonumber \\
\end{eqnarray}
where $e_{\Delta t}^{n+1}$ is the truncation error, $r=(t_{n+1}-s)$ and $b_j=(j+1)^{2-\gamma}-j^{2-\gamma}$. The reader may verify that
\begin{itemize}
  \item $b_j >0, j=0,1,2,\cdots,n$,
  \item $1=b_0>b_1>b_2>\cdots >b_n$ and $b_n\rightarrow 0$ as $n \rightarrow \infty$,
  \item $\sum\limits_{j=0}^n (b_j-b_{j+1})=(1-b_1)+\sum\limits_{j=1}^{n-1} (b_j-b_{j+1})+b_n=1$.
\end{itemize}
Substituting (\ref{2.2}) into (\ref{1.1}), we obtain the following temporal discretization
\begin{equation}\label{2.3}
\frac{1}{\Gamma(3-\gamma)}\sum\limits_{j=0}^n b_j \frac{u(x,t_{n+1-j})-2u(x,t_{n-j})+u(x,t_{n-1-j})}{\Delta t^\gamma}+\alpha u(x,t^{n+1})= \frac{\partial^2 u(x,t^{n+1})}{\partial x^2}+ f(x,t^{n+1}).
\end{equation}
Letting $\alpha_0=\frac{1}{(\Delta t)^{\gamma}\Gamma(3-\gamma)}, u^{n+1}=u(x,t_{n+1})$,  the last equation can be rewritten as
\begin{equation}\label{2.4}
\alpha_0(u^{n+1}-2u^{n}+u^{n-1})+\alpha_0 \sum\limits_{j=0}^n b_j ({u^{n+1-j}-2u^{n-j}+u^{n-1-j}})+\alpha u^{n+1}= \frac{\partial^2 u^{n+1}}{\partial x^2}+ f(x,t^{n+1}),
\end{equation}
where $n=0,1,\cdots M$.  It is observed that the term $u^{-1}$ will appear when $n=0$ or $j=n$. To eliminate $u^{-1}$, we utilize the initial condition to obtain
\begin{equation}\label{2.5}
u_t^0=\frac{u^1-u^{-1}}{2 \Delta t}.
\end{equation}
It follows then that $u^{-1}=u^1-2 \Delta t u_t^0$ or $u^{-1}=u^1-2 \Delta t \phi_2(x)$.
\section{Description of the Numerical Scheme}\label{s2}
In this section we derive the  cubic trigonometric B-spline collocation method (CuTBSM) for finding the numerical solution of time fractional diffusion-wave equation problem (\ref{1.1}). The solution domain $a\leq x\leq b$ is uniformly partitioned  by knots $x_{i}$ into $N$ subintervals $[x_{i}, x_{i+1}]$ of equal length $h=\frac{b-a}{N}$, $i=0,1,2,...,N-1$, where $a=x_{0}<x_{1}<...<x_{n-1}<x_{N}=b$. Our numerical approach for solving (\ref{1.1}) using trigonometric cubic B-splines is to seek an approximate solution $U(x,t)$ to the exact solution $u(x,t)$  in the following form  \cite{36,37}
\begin{equation}\label{3.1}
    U(x,t)=\sum\limits_{i=-1}^{N-1} c_i(t) TB_i^4 (x),
\end{equation}
where $c_i(t)$ are to be required for the approximate solution $U(x,t)$ to the exact solution $u(x,t)$.\\
The twice differentiable trigonometric basis functions $TB^4_i (x)$  \cite{abbas}  at the knots  $x_i$ are given by
\begin{equation}\label{3.2}
    TB^4_i (x)=\frac{1}{w}
    \begin{cases}
    p^3(x_i) & x\in [x_i,x_{i+1}]\\
    p(x_i)(p(x_i)q(x_{i+2})+q(x_{i+3})p(x_{i+1}))+q(x_{i+4})p^2(x_{i+1}),  & x \in [x_{i+1},x_{i+2}]\\
    q(x_{i+4})(p(x_{i+1})q(x_{i+3})+q(x_{i+4})p(x_{i+2}))+p(x_i)q^2(x_{i+3}), & x \in [x_{i+2},x_{i+3}]\\
    q^3(x_{i+4}), & x \in [x_{i+3},x_{i+4}]
    \end{cases}
\end{equation}
where\\
$\displaystyle{p\(x_i\)=\sin\(\frac{x-x_i}{2}\), q\(x_i\)=\sin\(\frac{x_i-x}{2}\), w=\sin\(\frac{h}{2}\) \sin\(h\) \sin\(\frac{3 h}{2}\).}$\\
Since there are three non zero terms at each knot notably $TB^4_{j-1}(x), TB^4_{j}(x)$ and $TB^4_{j+1}(x)$, therefore the approximation $u_j^n$ at the grid point $(x_j,t_n)$ to the exact solution at $n$ th time level is given as:
\begin{equation}\label{3.3}
    u(x_{j},t_{n})=u_j^n=\sum\limits_{j=i-1}^{i+1} c_j^n(t) TB^4_j (x).
\end{equation}
The  time dependent unknowns $c_j^n(t)$ are to be determined by making use of the initial and boundary conditions, and the collocation conditions on $TB^4_i (x)$. As a result we obtain the approximations $u_j^n$ together with its necessary derivatives as given below:
\begin{equation}\label{3.4}
    \begin{cases}
   \displaystyle{ u_j^n=a_1 c_{j-1}^n+a_2 c_{j}^n+ a_1 c_{j+1}^n}\\
    \displaystyle{ (u_x)_j^n=-a_3 c_{j-1}^n+a_3 \sigma_{j+1}^n}\\
    \displaystyle{ (u_{xx})_j^n=a_4 c_{j-1}^n+ a_5 c_{j}^n +a_4 c_{j+1}^n,}
    \end{cases}
\end{equation}
where\\
$\displaystyle{a_1=\csc\(h\) \csc\(\frac{3h}{2}\)\sin^2\(\frac{h}{2}\),}$\\
$\displaystyle{a_2=\frac{2}{1+2 \cos\(h\)},}$\\
$\displaystyle{a_3=\frac{3}{4} \csc \(\frac{3h}{2}\),}$\\
$\displaystyle{a_4=\frac{3+9\cos\(h\)}{4 \cos\(\frac{h}{2}\)-4 \cos\(\frac{5 h}{2}\)},}$\\
$\displaystyle{a_5=-\frac{3 \cot^2\(\frac{h}{2}\)}{2+4 \cos\(h\)}.}$\\\\
To obtain full discretization which relates the successive time levels and the unknowns $c_j^{n+1}$, we plug in the approximations $u_j^n$ and its  derivatives (\ref{3.4}) into the equation (\ref{2.4}). After some simplifications we arrive at the  following recurrence relation:\\
$\displaystyle{((\alpha_0+\alpha)a_1-a_4)c_{j-1}^{n+1} +((\alpha_0+\alpha)a_2-a_5)c_{j}^{n+1}+((\alpha_0+\alpha)a_1-a_4)c_{j+1}^{n+1}}=$
\begin{equation}\label{3.5}
         \begin{split}
        2\alpha_0(a_1 c_{j-1}^{n}+a_2 c_{j}^{n}+a_1 c_{j+1}^{n})-\alpha_0(a_1 c_{j-1}^{n-1}+a_2 c_{j}^{n-1}+a_1 c_{j+1}^{n-1})-\sum\limits_{k=1}^{n} b_k[a_1\(c_{j-1}^{n+1-k}-2c_{j-1}^{n-k}+c_{j-1}^{n-1-k}\)\\{}+a_2\(c_{j}^{n+1-k}-2c_{j}^{n-k}+c_{j}^{n-1-k}\)+a_1\(c_{j+1}^{n+1-k}-2c_{j+1}^{n-k}+c_{j+1}^{n-1-k}\)]+\alpha_0 f(x_j,t^{n+1}).
      \end{split}
     \end{equation}
The system (\ref{3.5}) contains $(N+1)$ linear equations in $(N+3)$ unknowns. To obtain two additional equations, the boundary conditions (\ref{1.3}) are utilized to obtain a unique solution of the problem. Consequently,  a  matrix system of dimension $(N + 3)\times(N + 3)$  is obtained which is a tri-diagonal  system.  Thomas Algorithm \cite{38} is then used to uniquely solve this system.
\subsection{Initial Vector $c^0$}
In order the commence the iteration process, it is required to find the initial solution vector $\displaystyle{c^0=[c_{-1}^0,c_0^0, \cdots, c_{N+1}^0]^T}$. The process of finding the initial vector involves the computation of  initial condition and its derivatives at the two boundaries as explained below \cite{abbas}:
\begin{enumerate}[(i)]
  \item $(u_j^0)_x=\frac{d}{dx}\phi_1(x_j),~~j=0$
  \item $u_j^0=\phi_1(x_j),~~ j=0,1,\cdots N$
  \item $(u_j^0)_x=\frac{d}{dx}\phi_1(x_j),~~ j=N$.
  \end{enumerate}
The above tri-diagonal system consists of $(N+3)$ linear equations in $(N+3)$ unknowns whose matrix form is given as:
\begin{equation}\label{3.6}
    A c^0=b,
\end{equation}
where
\begin{equation*}\label{21}
    A=\begin{bmatrix}
-a_3  & 0  & a_3 & \cdots & \cdots & \cdots & \cdots & 0 \\
a_1  & a_2  & a_1  & \ddots & && & \vdots \\
0 & a_{1}  & a_{2} & a_{1}  & \ddots & &  & \vdots \\
\vdots & \ddots & \ddots & \ddots & \ddots & \ddots &  & \vdots \\
\vdots & & \ddots & \ddots & \ddots & \ddots & \ddots& \vdots\\
\vdots  &  & & & \ddots & a_{1}  & a_{2}  &  a_{1}\\
0 & \cdots &  \cdots & \cdots & \cdots & -a_3 & 0 & a_3  \\
\end{bmatrix}
\end{equation*}
and $b=[\phi_1'(x_0), \phi_1(x_0), \cdots, \phi_1(x_N), \phi_1'(x_N)]^T.$
\section{Stability Analysis}\label{s4}
By  Duhamel's principle \cite{strik}, it follows that the solution to an inhomogeneous problem is the superposition of the solutions to homogeneous problems. As a consequence, a scheme is stable for inhomogeneous problem if it is stable for the homogeneous one.  It is sufficient to present the stability analysis for the scheme (\ref{3.5}) for the force free case ($f=0$) only.
The growth factor of  a Fourier mode is assumed to be $\rho_j^n$ and let $\tilde{\rho}_j^n$ be its approximation. Define $E_j^n=\rho_j^n-\tilde{\rho}_j^n$ which on substitution in  (\ref{3.5}), gives the following roundoff error equation:\\
$\displaystyle{((\alpha_0+\alpha)a_1-a_4)E_{j-1}^{n+1} +((\alpha_0+\alpha)a_2-a_5)E_{j}^{n+1}+((\alpha_0+\alpha)a_1-a_4)E_{j+1}^{n+1}}=$
\begin{equation}\label{4.1}
         \begin{split}
        2\alpha_0(a_1 E_{j-1}^{n}+a_2 E_{j}^{n}+a_1 E_{j+1}^{n})-\alpha_0(a_1 E_{j-1}^{n-1}+a_2 E_{j}^{n-1}+a_1 E_{j+1}^{n-1})-\\{}\sum\limits_{k=1}^{n} b_k[a_1\(E_{j-1}^{n+1-k}-2E_{j-1}^{n-k}+E_{j-1}^{n-1-k}\)+a_2\(E_{j}^{n+1-k}-2E_{j}^{n-k}+E_{j}^{n-1-k}\)\\{}+a_1\(E_{j+1}^{n+1-k}-2E_{j+1}^{n-k}+E_{j+1}^{n-1-k}\)].
      \end{split}
     \end{equation}
The error equation satisfies the boundary conditions
\begin{equation}\label{4.2}
E_0^k=\psi_1(t_k),E_N^k=\psi_2(t_k), k=0,1,\cdots,M
\end{equation}
and the initial conditions
\begin{equation}\label{4.3}
E_j^0=\phi_1(x_j),~~ (E_t)_j^0=\phi_2(x_j),~~  j=1,2,\cdots N.
\end{equation}
Define the grid function
$$E^k(x)=\begin{cases}
      E_j^k & x_j-\frac{h}{2} < x \leq x_j+\frac{h}{2}, j=1,\cdots,N-1 \\
      0 & a<x\leq \frac{h}{2} \text~{or}~ b-\frac{h}{2}<x\leq b.
   \end{cases}$$
Note that the Fourier expansion of $E^K(x)$ is
$$E^k(x)=\sum\limits_{m=-\infty}^\infty a_k(m)e^{\frac{i 2 \pi m x}{(b-a)}}, k=0,1,\cdots,M$$
where $a^k(m)=\frac{1}{(b-a)}\int\limits_a^b E^k(x)e^{\frac{-i 2 \pi m x}{(b-a)}} dx$. Let $$E^k=[E_1^k,E_2^k,\cdots,E_{N-1}^k]^T$$ and introduce the norm: $$\|E^k\|_2=\left(\sum\limits_{j=1}^{N-1} h |E_j^k|^2\right)^{\frac{1}{2}}=\left[\int\limits_a^b |E^k(x)|^2 dx\right]^{\frac{1}{2}}.$$
By  Parseval equality, it is observed that
$$\int\limits_a^b |E^k(x)|^2 dx=\sum\limits_{m=-\infty}^{\infty} |a_k(m)|^2,$$ so that the following relation is obtained
\begin{equation}\label{4.4}
\|E^k\|_2^2=\sum\limits_{m=-\infty}^\infty |a_k(m)|^2.
\end{equation}
Suppose that the equations (\ref{4.1}-\ref{4.3}) have solution of the form $E_j^n=\xi_n e^{i \beta j h}$, where  $i=\sqrt{-1}$ and $\beta$ is real. Substituting this expression in equation (\ref{4.1}) and dividing by $e^{i \beta j h}$ we obtain\\
$\displaystyle{((\alpha_0+\alpha)a_1-a_4)\xi_{n+1}e^{-i\beta h} +((\alpha_0+\alpha)a_2-a_5)\xi_{n+1}+((\alpha_0+\alpha)a_1-a_4)\xi_{n+1}e^{i\beta h}=}$
\begin{equation}\label{4.5}
         \begin{split}
        2\alpha_0(a_1 \xi_{n}e^{-i\beta h}+a_2 \xi_{n}+a_1\xi_{n}e^{i\beta h})-\alpha_0(a_1 \xi_{n-1}e^{-i\beta h}+a_2 \xi_{n-1}+a_1\xi_{n-1}e^{i\beta h})-\\{}\sum\limits_{k=1}^{n} b_k[a_1\(\xi_{n+1-k}e^{-i\beta h}-2\xi_{n-k}e^{-i\beta h}+\xi_{n-1-k}e^{-i\beta h}\)+a_2\(\xi_{n+1-k}-2\xi_{n-k}+\xi_{n-1-k}\)\\{}+a_1\(\xi_{n+1-k}e^{i\beta h}-2\xi_{n-k}e^{i\beta h}+\xi_{n-1-k}e^{i\beta h}\)].
      \end{split}
     \end{equation}
Using the relation $e^{-i \beta h}+e^{i \beta h}=2\cos(\beta h)$ and grouping like terms, we obtain the following relation
\begin{equation}\label{4.6}
\xi_{n+1}=\frac{2}{\nu} \xi_{n}-\frac{1}{\nu} \xi_{n-1}-\frac{1}{\nu}\sum\limits_{k=1}^{n} b_k\(\xi_{n+1-k}-2\xi_{n-k}+\xi_{n-1-k}\),
\end{equation}
where $\nu=1+\frac{2(a_1\alpha-a_4)\cos(\beta h))+(a_2\alpha-a_5)}{\alpha_0(2 a_1 \cos(\beta h)+a_5)}$. Obviously $\nu \geq 1$.
\begin{prop}\label{aa}
If $\xi_n$ is the solution of equation (\ref{4.6}), then $|\xi_n|\leq 2|\xi_0|$, $n=0,1,\cdots, T \times M$.
\end{prop}
\noindent
\textbf{Proof}. Mathematical induction is used to prove the result. For $n=0$, we have from equation $(\ref{4.6})$ that $\xi_1=\frac{2}{\nu} \xi_0$. Since $\nu \geq 1$, we have
$$|\xi_1|=\frac{2}{\nu} |\xi_0|\leq 2|\xi_0|.$$
Now suppose that $|\xi_n|\leq 2|\xi_0|$, $n=1,\cdots,T \times M-1 $ so that from (\ref{4.6}), we obtain
\begin{eqnarray}
  |\xi_{k+1}|&\leq& 2\frac{|\xi_n|}{\nu}-\frac{|\xi_{n-1}|}{\nu}-\frac{1}{\nu}\sum\limits_{k=1}^{n}b_k\(|\xi_{n+1-k}|-2|\xi_{n-k}|+|\xi_{n-1-k}|\)\\ \nonumber
    &\leq& 4\frac{|\xi_0|}{\nu}-2\frac{|\xi_{0}|}{\nu}-\frac{2}{\nu}\sum\limits_{k=1}^{n}b_k\(|\xi_{0}|-2|\xi_{0}|+|\xi_{0}|\)\\ \nonumber
        &\leq& 2|\xi_0|.
\end{eqnarray}
This completes the proof.
\begin{theorem}
The collocation scheme (\ref{3.5}) is unconditionally stable.
\end{theorem}
\noindent
\textbf{Proof}. Using formula (\ref{4.4}) and Proposition \ref{aa}, we  obtain $$\|E^k\|_2\leq 2\|E^0\|_2, k=0,1,\cdots M $$  which establishes that the  scheme is unconditionally stable.
\section{Numerical Results and Discussion}\label{s5}
In this section, numerical experiments are carried out for the time fractional diffusion-wave equation (\ref{1.1}) with initial (\ref{1.2}) and boundary conditions (\ref{1.3}). To check the efficiency and accuracy of the method is checked by calculating the error norms $L_2$ and $L_\infty$  given by
$$L_2=\|U^{\text{exact}}-U_{N}\|_2\simeq \sqrt{h\sum\limits_{j=-1}^{N+1}|U_j^{\text{exact}}-(U_N)_j|}$$
and $$L_\infty=\|U^{\text{exact}}-U_{N}\|_\infty \simeq \max\limits_j |U_j^{\text{exact}}-(U_N)_j||$$\\
respectively. We compare the numerical solutions obtained by CuTBSM for one dimensional fractional diffusion equations (\ref{1.1}) with known exact solutions. Numerical calculations are carried out by using Mathematica 9 on an Intel$\circledR$Core$^{TM}$ i5-2410M CPU $@$ 2.30 GHz with 8GB RAM and 64-bit operating system (Windows 7).
\begin{exmp}
As a first experiment we consider the following fractional diffusion-wave equation for $\alpha=0$
\begin{eqnarray}\label{5.1a}
  \frac{\partial^ \gamma }{\partial t^\gamma} u(x,t) &=&\frac{\partial ^2}{\partial x^2} u(x,t)+ f(x,t),~~\qquad\qquad (x,t)\in [0,1]\times [0,T] \nonumber \\
    u(x,0) &=& 0,~ u_t(x,0)=-\sin(\pi x) \nonumber\\
  u(0,t) &=&u(1,t)=0,
\end{eqnarray}
\end{exmp}
where the source term is $ f(x,t)=\frac{2 t^{2-\gamma}\sin(\pi x)}{\Gamma[3-\gamma]}+(t^2-t)\sin(\pi x)\pi^2$. The exact solution of the problem is $u(x,t)=\sin(\pi x)(t^2-t)$ \cite{khader}.
Figure 1 compares the graphs of the exact and approximate solutions with different values of $\gamma, h$ and $\Delta t$ at different time levels. In Figure 2, we exhibit the absolute error profiles at different time levels. Figure 3 compares the graphs of exact and approximate solutions using our scheme with those obtained in \cite{khader} at time $t=2$. It is observed that our scheme gives much better accuracy. Figure 4 shows 3D plots of exact and approximates solutions at time $t=0.1$. In Tables 1-2, the maximum errors obtained are compared with those of Hermite Formula (HF) \cite{khader} for different values of $\gamma$ to demonstrate that our scheme is more accurate and gives accuracy of $10^{-7}$. In tables 3-4, error norms are computed for various values of parameters to further confirm the accuracy and efficiency of the presented scheme.
\FloatBarrier
\begin{figure}[h!]
\begin{center}
%\subfigure{\includegraphics[width=2.9in]{fig1d1}} \hspace{0.3in}
%\subfigure{\includegraphics[width=2.9in]{fig1c1}}
\caption{Comparison of exact and numerical solutions for Example 1 when $h=\frac{1}{80}, \Delta t=\frac{1}{100}$ and $\gamma=1.75$.}
\end{center}
\end{figure}
\FloatBarrier
\begin{figure}[h!]
\centering
\caption{Error profiles for $\gamma=1.5, h=\frac{1}{150},\Delta t=\frac{1}{120}$, at different time levels for Example 1.}
\end{figure}
\FloatBarrier
\begin{figure}[h!]
\begin{center}
%\subfigure[Behavior of numerical solution obtained in \cite{khader}]{\includegraphics[width=2.9in]{fig6a1}} \hspace{0.3in}
%\subfigure[Comparison of exact(solid line) and approximate solution(dots) by present scheme]{\includegraphics[width=2.9in]{fig61}}
\caption{Behavior of numerical solutions for for Example 1 when $\gamma=1.5, h=\frac{1}{150}, \Delta t=\frac{1}{120}$ and $t=2$.}
\end{center}
\end{figure}
\FloatBarrier
\begin{figure}[h!]
\begin{center}
%\subfigure{\includegraphics[width=2.9in]{fig5b1}} \hspace{0.3in}
%\subfigure{\includegraphics[width=2.9in]{fig5a1}}
\caption{3D comparison of exact (left) and numerical solutions (right) for for Example 1 when $h=\frac{1}{64}, \Delta t=\frac{1}{100}$ and $\gamma=1.5$ at $t=1$.}
\end{center}
\end{figure}
\FloatBarrier
\begin{exmp}\label{exp1}
As a second experiment consider the time fractional diffusion-wave equation
\begin{eqnarray}\label{5.1}
  \frac{\partial^ \gamma }{\partial t^\gamma} u(x,t)+u(x,t)  &=&\frac{\partial ^2}{\partial x^2} u(x,t)+ f(x,t),~~\qquad\qquad (x,t)\in [0,1]\times [0,T] \nonumber \\
    u(x,0) &=& 0,~ u_t(x,0)=0 \nonumber\\
  u(0,t) &=&u(1,t)=0,
\end{eqnarray}
where the forcing term $f(x,t)$ is supposed to be $\displaystyle{f(x,t)=\frac{2t^{2-\gamma}x(1-x)}{\Gamma[3-\gamma]}+t^2x(1-x)+2 t^2}$. The exact solution of the  problem (\ref{5.1}) is $u(x,t)=t^2 x (1-x)$.
\end{exmp}
The proposed scheme is applied to solve this problem. Figure 5 shows the graphs of the exact and approximate solutions at different time levels with $\gamma=1.5$, $\Delta t=0.01$ and $t=1$ for  $N=80$. The absolute error profile is plotted at different time levels in Figure 6 (with $N=60, \Delta t=0.001$). The approximate and exact solutions are shown in Figure 7 by fixing values of different parameters. The absolute errors at different points $(x_j,t_j)\in [0,1]\times [0,1]$ for different values of $\gamma$ chosen in the range $1<\gamma \leq 2$ are tabulated in Table \ref{T1}. From Figure 7 and Table \ref{T1}, it is clear that the proposed scheme is very accurate and efficient. It is worthwhile to note that the numerical solutions are in excellent agreement with the exact solutions for many values of $\gamma$.
\begin{figure}[h!]
\centering
\caption{ The comparison of exact (lines) and approximate solutions for Example 2 when $\gamma=0.5$, at $t=0.5$ (rectangles), $t=0.75$ (circles) and $t=1$ (stars) with $N=80,~\Delta t=0.01$.}
\end{figure}
\FloatBarrier
\begin{figure}[h!]
\centering
\caption{Error profiles for Example 2 at different time levels when $N=60, \Delta t=0.01$.}
\end{figure}
\FloatBarrier
%\begin{table}[h!]
%\setlength{\tabcolsep}{8pt}
%\centering \caption{Absolute errors of Example 2 for many values of $\gamma$ at different points}\label{T1}
%\begingroup\setlength{\fboxsep}{0pt}
%\vspace{2mm}
%\colorbox{lightgray}{%
%\begin{tabular}{c c ccc c c}
%\hline\hline $(x,t)$  & $\gamma=1.1$ & $\gamma=1.3$ & $\gamma=1.5$ & $\gamma=1.7$ & $\gamma=1.9$\\
%\hline
%(0.1,0.1)  & 9.5133e-09 & 6.6004e-09 & 4.4920e-09 & 2.9885e-09 & 1.9326e-09\\
%(0.2,0.2)  & 1.0530e-07 & 7.9127e-08 & 5.7844e-08 & 4.1404e-08 & 2.8903e-08\\
%(0.3,0.3)  & 9.6665e-07 & 3.3461e-07 & 2.5678e-07 & 1.9243e-07 & 1.4105e-07\\
%(0.4,0.4)  & 1.0813e-06 & 9.1574e-07 & 7.3594e-07 & 5.7117e-07 & 4.3402e-07\\
%(0.5,0.5)  & 2.2190e-06 & 1.6516e-06 & 1.6516e-06 & 2.2190e-06 & 1.0367e-06\\
%(0.6,0.6)  & 3.9341e-06 & 3.1570e-06 & 3.1570e-06 & 3.9341e-06 & 2.1091e-06\\
%(0.7,0.7)  & 6.3114e-06 & 5.3801e-06 & 5.3809e-06 & 6.3114e-06 & 3.8408e-06\\
%(0.8,0.8)  & 9.4176e-06 & 8.4176e-06 & 8.4176e-06 & 9.4176e-06 & 6.4421e-06\\
%(0.9,0.9)  & 1.3302e-05 & 1.2877e-05 & 1.2324e-05 & 1.3302e-05 & 1.3302e-05\\
% \hline\hline
%\end{tabular}%
%}\endgroup
%\end{table}
\FloatBarrier
\begin{figure}[h!]
\begin{center}
%\subfigure{\includegraphics[width=2.9in]{fig1c}} \hspace{0.3in}
%\subfigure{\includegraphics[width=2.9in]{fig1d}}
\caption{ 3D space-time graph of exact and numerical solution for Example 2}
\end{center}
\end{figure}
\begin{exmp}
As a last example consider the  time fractional diffusion-wave equation
\begin{eqnarray}\label{5.2}
  \frac{\partial^ \gamma }{\partial t^\gamma} u(x,t)+u(x,t)  &=&\frac{\partial ^2}{\partial x^2} u(x,t)+ f(x,t),~~\qquad\qquad (x,t)\in [0,1]\times [0,T] \nonumber \\
    u(x,0) &=& 0,~ u_t(x,0)=0 \nonumber\\
  u(0,t) &=&0, u(1,t)=t^2 \sinh(1),
\end{eqnarray}
where the source term is  $\displaystyle{f(x,t)=\pi \frac{2\sinh(x)t^{2-\gamma}}{\Gamma[3-\gamma]}}+(1-\pi)t^2\sinh(x)$. The exact solution of the problem (\ref{5.2}) is $\displaystyle{u(x,t)=t^2 \sinh(x).}$
\end{exmp}
The above problem is solved by using the proposed scheme. Figure 8 exhibits the graphs of the exact and approximate solutions at different time levels with $\gamma=1.5$, $\Delta t=0.01$ and $t=1$ for  $N=80$.  By taking  $N=40, \Delta t=0.001$, the absolute error profile is plotted at different time levels in Figure 9. The 3D exact and numerical solutions are shown in Figure 10 by fixing values of different parameters. In table 6, we present the absolute errors at different points $(x_j,t_j)\in [0,1]\times [0,1]$ for different values of $\gamma$ chosen in the range $1<\gamma \leq 2$. A comparison between the obtained results and those of radial basis functions (RBF) \cite{vahid} is given in Tables 7-9. From Figure 10 and Tables 6-9, it is clear that the proposed scheme is very accurate and efficient. It is noticed that the numerical solutions are in close agreement  with the exact solutions for many values of $\gamma$.
\begin{figure}[h!]
\centering
\caption{ The comparison of exact (lines) and approximate solutions for Example 3 when $\gamma=0.5$, at $t=0.5$ (rectangles), $t=0.75$ (circles) and $t=1$ (stars) with $N=80,~\Delta t=0.01$.}
\end{figure}
\FloatBarrier
\begin{figure}[h!]
\centering
\caption{Error profiles for Example 3 at different time levels when $N=60, \Delta t=0.01$.}
\end{figure}
\FloatBarrier
\begin{figure}[h!]
\begin{center}
%\subfigure{\includegraphics[width=2.9in]{fig2c}} \hspace{0.3in}
%\subfigure{\includegraphics[width=2.9in]{fig2d}}
\caption{3D space-time graph of exact and numerical solution for Example 3}
\end{center}
\end{figure}
%\begin{table}[h!]
%\setlength{\tabcolsep}{8pt}
%\centering \caption{Absolute errors of Example 3 for several values of $\gamma$ at different points} \label{TT}
%\begingroup\setlength{\fboxsep}{0pt}
%\vspace{2mm}
%\colorbox{lightgray}{%
%\begin{tabular}{c c ccc c c}
%\hline\hline $(x,t)$  & $\gamma=1.1$ & $\gamma=1.3$ & $\gamma=1.5$ & $\gamma=1.7$ & $\gamma=1.9$\\
%\hline
%(0.1,0.1)  & 1.2498e-09 & 6.0892e-10 & 2.9307e-10 & 1.3992e-10 & 6.6457e-11\\
%(0.2,0.2)  & 2.0339e-08 & 1.1879e-08 & 6.5709e-09 & 3.5960e-09 & 1.9547e-09\\
%(0.3,0.3)  & 9.7339e-08 & 6.6352e-08 & 4.0629e-08 & 2.4101e-08 & 1.4186e-08\\
%(0.4,0.4)  & 2.8201e-07 & 2.1618e-07 & 1.4732e-07 & 9.3061e-08 & 1.7977e-08\\
%(0.5,0.5)  & 6.2459e-07 & 5.1997e-07 & 3.9209e-07 & 2.6549e-07 & 1.7292e-07\\
%(0.6,0.6)  & 1.1728e-06 & 1.0335e-06 & 8.4665e-07 & 6.2213e-07 & 4.2256e-07\\
%(0.7,0.7)  & 1.9709e-06 & 1.8065e-06 & 1.5758e-06 & 1.2559e-06 & 8.9991e-07\\
%(0.8,0.8)  & 3.0599e-06 & 2.8817e-06 & 2.6324e-06 & 2.2492e-06 & 1.7266e-06\\
%(0.9,0.9)  & 4.4761e-06 & 4.4761e-09 & 4.0559e-06 & 3.6618e-06 & 3.0216e-06\\
% \hline\hline
%\end{tabular}%
%}\endgroup
%\end{table}
\FloatBarrier
\FloatBarrier
\section{Concluding Remarks}\label{s6}
This study presents a finite difference scheme with a combination of  cubic trigonometric B-spline basis for the time fractional fractional diffusion-wave equation with reaction term. This algorithm is based on a discretization using  finite difference formulation  for the Caputo sense. The cubic trigonometric B-spline basis functions have been used to approximate derivatives in space. The scheme provides accuracy of $10^{-11}$ and the obtained numerical results are in superconformity with the exact solutions. A special attention has been given to study the stability of the scheme by using a procedure similar to Von Neumann stability analysis. The scheme is shown to be unconditionally stable, whereas the scheme of \cite{khader} is conditionally stable.

%%%%%%%%%%%%%%%

\begin{thebibliography}{33}
\bibitem{pod} I. Podlubny, \emph{Fractional Differential Equations}, Academic Press, 1999.
\bibitem{mai}
F. Mainardi, \emph{Fractals and Fractional Calculus Continum Mechanics}, Springer Verlag, (1997), 291-348.
\bibitem{hil}
R. Hilfer, \emph{Applications of Fractional Calculus in Physics}, World Scientific, Singapore, 2000.
\bibitem{kilbas}
A.A. Kilbas, H.M. Srivastava, and J.J. Trujillo, \emph{Theory and Applications of Fractional Differential Equations}, Elsevier Science B.V., Amsterdam, 2006.
\bibitem{sok}
I. M. Sokolov, J. Klafter, and A. Blumen, \emph{Fractional kinetics}, Phys. Today, 55 (2002), 48-54.
\bibitem{die}
K. Diethelm, A. D. Freed,\emph{ On Solution of nonlinear fractional order differential equations used in modelling of viscoplasticity in: Scientific Computing in Chemical Engineering II: Computational Fluid Dynamics, Reaction Engineering and Molecular Properties}, Springer Verlag, Heidelberg, 1999, 217-224.
\bibitem{bek}
P. Beker-Keren,  M. M. Meerschaert, H. P. Scheffler, \emph{Limit theorem for continuous-time random walks with two time scales}, J. Appl. Prob. 41 (2004) 455-466.
\bibitem{meer}
M. M. Meerschaert, Y. Zhang, B. Baeumerc \emph{Particle tracking for fractional diffusion with two time scales}, Comput. Math. Appl. 59 (2010) 1078-1086.
\bibitem{iom}
A. Iomin, S. Dorfman, and L. Dorfman, \emph{On tumor development: fractional transport approach}, http://arxiv.org/abs/qbio/0406001.
\bibitem{scale1}
R.Gorenflo, F. Mainradi et al \emph{Fractional calculus and continuous-time finance. III, The diffusion limit.},Mathematical Finance, Trends in Math. (2001) 171-180.
\bibitem{scale2}
M. M. Meerschaert, E. Scalas, F. Mainradi \emph{Coupled continuous time random walks in finance}, Physica A, 370 (2006), 114-118.
\bibitem{scale3}
M. Raberto, E. Scalas, F. Mainradi \emph{Waiting-times and returns in high-frequency financial data: an empirical study}, Physica A, 314 (2002), 749-755.
\bibitem{zeng}
F. Zeng \emph{Second-Order stable finite difference schemes for the time-fractional diffusion-wave equation}, J. Sci. Comput., (2014), 1-20.
\bibitem{khader}
M. M. Khader, M. H. Adel \emph{Numerical solutions of fractional wave equations using an efficient class of FDM based on the Hermite formula}, Advances in Difference Equations, (2016), 1-10.
%\bibitem{khader2}
%N H. Sweilam, M. M. Khader, M. H. Adel \emph{On the stability analysis of weighted average finite difference method for fractional wave equation},Fract. Differ. Calc., 2(1) (2012), 17-29.
\bibitem{vahid}
Z. Avazzadeh, V. R. Hosseini, W..Chen  \emph{Radial basis functions and FDM for solving fractional diffusion-wave equation},IJST, 38(A3) (2014), 205-212.
\bibitem{psku}
A. V. Pskhu  \emph{The fund amental solution of a diffusion-wave e quation of fractional order}, Izvestiya: Mathematics, 73(2) (2009), 351-392.
\bibitem{pov}
Y. Povstenko \emph{Neumann boundary-value problems for a time-fractional diffusion-wave equation in a half plane}, camwa, 64 (2012),  3183–3192.
%\bibitem{william}
%W. McLean, K. Mustapha \emph{A second order accurate numerical method for a fractional wave equation}, Numer. Math., 105 (2007), 481–510.
\bibitem{lie}
A. Liemert, A. Kienle \emph{Time-fractional wave-diffusion equation in an inhomogeneous half-space}, J. Phys. A: Math. Theor. , 48 (2015), 1–19.
%\bibitem{khader3}
%N. H. Sweilam, M. M. Khader, H. M. Almarwm \emph{Numerical studies for the variable order nonlinear fractional wave equation},Fract.  Calc. Appl. Anal., 15(4) (2012), 669-683.
%\bibitem{haddar}
% H. Haddar, J. -R. Li,, D. Matignon\emph{Efficient solution of a wave equation with fractional-order dissipative terms},cam, 234 (2010), 2003-2010.
\bibitem{ren}
J. Ren, Z. Sun\emph{Efficient numerical solution of the multi-term time fractional diffusion-wave equation },eajam, 5 (1) (2015), 1-28.
\bibitem{jin}
B. Jin, R. Lazarov, Y. Liu, Z.Zhou\emph{The Galerkin finite element method for a multi-term time-fractional diffusion equation },J. Comp. Phy., 281 (2015), 825-843.
\bibitem{36}
P.M. Prenter, \emph{Splines and Variational Methods}, John Wiley and Sons, 1989.
\bibitem{37}
C. de Boor, \emph{A Practical Guide to Splines}, Springer-Verlag, 1978.
\bibitem{abbas}
M. Abbas, A. A. Majid, A. I. M. Ismail, and A. Rashid, \emph{The application of cubic trigonometric B-spline to the numerical solution of the hyperbolic problems}, Appl. Math. Comput. 239 (2014), 74-88.
\bibitem{38}
R.L. Burdern, J.D. Faires, \emph{Numerical Analysis}, eighth ed., Brooks Cole, 2004.
\bibitem{strik}
J. C. Strikwerda, \emph{Finite Diffrence Schemes and Partial Diffrential Equations}, Society for Industrial and Applied Mathematics,Philadelphia, Pa, USA, 2nd edition, 2004.
\end{thebibliography}
\end{document}